\documentclass[11pt,a4paper]{amsart}
\title[Gröbner walk in OSCAR]{The \texttt{GroebnerWalk.jl} package for \texttt{OSCAR}}
\author{Kamillo Ferry}
\address{Technical University Berlin, Germany}
\email{ferry@math.tu-berlin.de}

\author{Francesco Nowell}
\address{Technical University Berlin, Germany}
\email{nowell@math.tu-berlin.de}


\usepackage{macros}
\usepackage{graphicx} 

\usepackage{siunitx}

\usepackage{standalone}
\usepackage[frozencache=true,cachedir=minted-cache]{minted}
\usepackage{listings}
\usepackage{amsmath}
\usepackage{amssymb}
\usepackage{amsthm}
\usepackage{amsxtra}
\usepackage{mathtools}

\usepackage{enumitem}

\usepackage{xcolor}
\usepackage{tikz}
\usetikzlibrary{arrows,decorations.pathreplacing,shapes}
\usetikzlibrary{positioning,fit}
\usetikzlibrary{graphs,graphs.standard,calc}
\usepackage{pgfplots}
\pgfplotsset{compat=1.18}

\usepackage{eqnarray}

\usepackage[mathscr]{euscript}
\usepackage{stmaryrd}

\tikzset{cross/.style={cross out, draw,%
         minimum size=2*(#1-\pgflinewidth),%
         inner sep=0pt, outer sep=0pt}}

\definecolor{bluevar}{rgb} {0.2,0.37,0.64}
\definecolor{redvar}{rgb}{0.64, 0, 0}
\definecolor{ltred}{rgb}{0.64, 0, 0}

\tikzset{
    circlenode/.style={circle, draw=black!85, fill=black!35, scale=.5},
    endpoint/.style={circlenode,redvar},
    intermediate/.style={circlenode},
    ray/.style = {thick, color=bluevar},
}

\lstdefinelanguage{Julia}%
  {morekeywords={abstract,break,case,catch,const,continue,do,else,elseif,%
      end,export,false,for,function,immutable,import,importall,if,in,%
      macro,module,otherwise,quote,return,switch,true,try,type,typealias,%
      using,while},%
   sensitive=true,%
   alsoother={\$},%
   morecomment=[l]\#,%
   morecomment=[n]{\#=}{=\#},%
   morestring=[s]{"}{"},%
   morestring=[m]{'}{'},%
}[keywords,comments,strings]

\usepackage{setspace}
\usepackage{csquotes}

\usepackage{hyperref}
\hypersetup{
  bookmarksopen,
  bookmarksnumbered,
  pdftitle={Gröbner walk in OSCAR}
  pdfauthor={Kamillo Ferry, Francesco Nowell},
  pdfkeywords={project,report},
  pdfstartview={FitH},
}

\usepackage[headings]{fullpage}

\usepackage[backend=biber,style=numeric,backref,backrefstyle=none]{biblatex}
\AtNextBibliography{\small}

\usepackage[nameinlink]{cleveref}

\theoremstyle{plain}
\newtheorem{theorem}{Theorem}

\newtheorem{lemma}[theorem]{Lemma}

\theoremstyle{definition}

\newtheorem{example}[theorem]{Example}

\theoremstyle{remark}

\newcommand{\citeauthorandcite}[1]{\citeauthor{#1}~\cite{#1}}

\newcommand{\R}{\mathbb{R}}

\newcommand{\Q}{\mathbb{Q}}

\newcommand{\tless}{<}

\newcommand{\supp}{\mathrm{supp}}

\newcommand{\lexord}{\tless_{\mathrm{lex}}}
\newcommand{\grevlexord}{\tless_{\mathrm{degrevlex}}}

\newcommand{\inless}{\initial_\tless}
\newcommand{\inw}{\initial_\omega}

\newcommand{\w}{\omega}

\begin{document}
\onehalfspace
\begin{abstract}
   Gröbner bases are a central tool in computational algebra, but it is well-known that their ease of computation rapidly deteriorates with increasing number of variables and/or degree of the input generators.
   Due to the connection between polyhedral geometry and Gröbner bases through the \emph{Gröbner fan},
   one can attempt an incremental approach to compute Gröbner bases. 
   First computing a Gröbner basis with respect to an 
   `easy' term ordering and transforming that result to a
   Gröbner basis with respect to the desired term ordering
   by using information about this polyhedral fan is done by 
   a family of algorithms termed as \emph{Gröbner walk}.
   We implemented two variants of the Gröbner walk
   in the computer algebra system \texttt{OSCAR} and compared their 
   performance with classical Gröbner basis methods already found in \texttt{OSCAR}.
\end{abstract}

\maketitle
\section{Introduction}
The \emph{Gröbner walk} is an approach to reduce the computational 
complexity of Gröbner basis computations first proposed by 
\citeauthor{Collart:GWalk:1997}~\cite{Collart:GWalk:1997}.
Algorithms of this type work by exploiting the geometry of the \emph{Gröbner fan} which has been introduced by Mora and Robbiano \cite{Mora.Robbiano:1998}. This is a polyhedral fan associated to a polynomial ideal $I$, the maximal cones of which are in one-to-one correspondence with the Gröbner bases of $I$ as one varies over all term orderings. These algorithms belong to the wider class of Gröbner basis methods which operate incrementally via subsequent \emph{conversions} (e.g FGLM \cite{Faugere:1993} or Hilbert-driven Buchberger \cite{Traverso:1997}).  

The aim of the Gröbner walk is to compute a Gröbner basis for an ideal \(I\) given two term orderings of a polynomial ring,
the start ordering \(<_s\) and the target ordering \(<_t\).
The algorithm starts by computing a \emph{Gröbner basis} with respect to \(<_s\) and the corresponding cone in the Gröbner fan. It then computes a boundary point in the direction of the cone corresponding to desired target ordering. 
This point $\w$ is then used to retrieve a basis of the ideal of initial forms $\inw(I)$, which in the Gröbner fan corresponds to the lower-dimensional cone on which $\w$ lies. This basis is subsequently lifted to the Gröbner basis of $I$ corresponding to the adjacent full-dimensional cone. This process is repeated until one obtains the target Gröbner basis after finitely many steps.
Our package \texttt{GroebnerWalk.jl} implements the original algorithm described in \cite{Collart:GWalk:1997} (which we refer to as the \emph{standard walk}) as well as the \emph{generic Gröbner walk} of \citeauthor{GenericWalk:2007}, which combines methods of the standard walk with techniques from symbolic computation. It adds to the array of Gröbner basis algorithms already accessible in the Julia ecosystem via the \texttt{msolve}/\texttt{AlgebraicSolving.jl}, \texttt{OSCAR} \cite{OSCAR,OSCAR-book:2025} and \texttt{Groebner.jl} packages.

\section{Basics of Gröbner fans}
A \emph{term ordering} on \(R = k[x_1,\dots,x_n]\) is a relation \(\tless\) on the monomials of \(R\) which is a strict total well-ordering 
that satisfies \[
  \text{for all } \alpha,\beta,\gamma\in\NN^n \ \colon \  x^\alpha \tless x^\beta \implies x^{\alpha+\gamma} \tless x^{\beta+\gamma}.
\] For any non-zero polynomial \(f\in R\) and term ordering \(\tless\), there is a unique maximal term \(c_\alpha x^\alpha\) with respect to
\(\tless\) which is called the \emph{leading term} \(\inless(f)\).

The \emph{initial ideal} of a non-zero ideal \(I \trianglelefteq R\) \wrt{} \(\tless\) is \[
  \inless(I) = \langle\, \inless(f) \mid f\in I \,\rangle.
\] Sometimes, this is also referred to as \emph{leading ideal} in literature and also in \texttt{OSCAR}. 
A finite set $G = \{ g_1, \dots g_s\} \subset R$ of polynomials is called a \emph{Gröbner basis} 
for an ideal $I$ w.r.t $<$ if \(I=\langle G\rangle\) and
\begin{equation*}
     \langle\, \inless(g_1) \dots \inless(g_s) \,\rangle = \inless(I).
\end{equation*}

\begin{example}\label{ex:running-example}
    We demonstrate the theory and the corresponding functions in \texttt{OSCAR} using the ideal
    \begin{equation}
        I = \langle y^4+ x^3-x^2+x,x^4\rangle \triangleleft \QQ[x,y]
    \end{equation} as running example.

    The following code snippet defines this ideal and calculates the initial ideal with respect to
    the lexicographic ordering \(\lexord\). 
\begin{minted}{jlcon}
julia> using Oscar
julia> R, (x,y) = QQ[:x, :y]; I = ideal([y^4+ x^3-x^2+x,x^4]);
julia> leading_ideal(I; ordering=lex(R))
Ideal generated by
  y^16
  x
\end{minted}
\end{example}
Given a vector $\w \in \R^n_{\geq 0}$ and a term ordering $<$, we can define the \emph{weight ordering} \[
  \alpha\tless_\w\beta
  \colon\iff 
  \langle\omega,\alpha\rangle < \langle\omega,\beta\rangle 
  \textit{ or } 
  (\langle\omega,\alpha\rangle = \langle\omega,\beta\rangle 
  \textit{ and }
  x^\alpha\tless x^\beta).
\]

This is also a term ordering, sometimes referred to as the \emph{refinement} of the weight vector $\w$ by $<$. Upon relaxing the comparison with $<$, we obtain the \emph{partial weight ordering} $\prec_\w$: 
\[ \alpha \prec_\w \beta : \iff \langle \alpha , \w \rangle < \langle \beta , \w \rangle , \]
which in general is not a term ordering. 
For example, in the setting of \Cref{ex:running-example} a refinement of the weight vector \((2,1)\) by \(\lexord\) 
may be obtained in \texttt{OSCAR} as follows.
\begin{minted}{jlcon}
julia> weight_ordering([2,1], lex(R))
matrix_ordering([x, y], [2 1])*lex([x, y])
\end{minted}
The \emph{initial form} \(\inw(f)\) of \(f\) \wrt{} \(\omega\) is defined as the sum of 
the terms of \(f\) which are maximal with respect to $\prec_\w$. The \emph{generalized initial ideal} of \(I\) \wrt{} \(\omega\) is the ideal generated by the initial forms: \[
  \inw(I) = \langle\, \inw(f) \mid f\in I \,\rangle.
\] Note that \(\inw(I)\) is not a monomial ideal in general. For example,
for the ideal \(I\) from \Cref{ex:running-example}, the weight vector \((4,3)\in\R^2\) yields a partial weight ordering 
\(\prec_\w\) such that \(\inw(I) = \langle\, x^3 + y^4, x^4 \,\rangle \) which is not monomial.

For a fixed ideal $I$ and term ordering $<$ we say that a weight vector $\omega$ \emph{represents} $<$ if $\initial_\w(I) = \initial_< (I)$. If a Gröbner basis $G$ w.r.t. $<$ is given, a necessary and sufficient condition for $\w$ to represent $<$ is that $\inw(g) = \initial_<(g) $ holds for all $ g \in G$. The crucial connection to polyhedral geometry is that the set of all weight vectors representing a given term ordering $<$  lie in the relative interior of a full-dimensional polyhedral cone in $\R^n$ with integer generators. Upon taking the closure of this cone, and then the union of all such cones varying over all term orderings, one obtains a full-dimensional rational polyhedral fan, called the \emph{Gröbner fan} of $I$. We denote this by $\mathbb{G}(I)$. Another key notion is that of a \emph{marked Gröbner basis}, which is a Gröbner basis $G = \{ g_1, \dots g_s\}$ for $I$ w.r.t. $<$ with the following additional properties: 
\begin{enumerate}[label = (\textit{\roman*})]
\item $G$ is minimal w.r.t inclusion, i.e. 
\begin{center}
     $ \inless(G \setminus \{g_i \}) \subsetneq \inless(I)\ \ $  for all $i \in \{ 1, ..., s \}$.
\end{center}
\item $G$ is monic, i.e. the coefficient of \(\inless(g_i)\) is equal to \(1\) for all $i \in \{ 1, ..., s \} $. 
\item $G$ is reduced, i.e. 
\begin{center}
    for all $i, j \in \{ 1, ..., s \}$ with $i \neq j$, no term of $g_i$ is divisible by $\inless(g_j)$.
\end{center}
\item The leading terms of elements of $G$ are \emph{marked}, in the sense that each $g \in G$ is formally encoded as the pair $(g, x^\alpha)$, where $x^\alpha = \inless(g)$.
\end{enumerate}
Marked Gröbner bases are unique in the sense that if two term orderings give rise to the same initial ideal, then their marked Gröbner bases coincide. 

Some key theoretical results behind the Gröbner walk are stated below. Detailed proofs and additional context may be found in Chapters 1 and 2 of \cite{Nowell:2025}.
\begin{theorem} \label{maintheorem}
For a non-zero ideal $I \triangleleft R$, the following sets are in one-to-one correspondence: 
\\
\begin{center}
$\Biggl\{$ \parbox[c][2.5em][c]{3.8cm}{\centering $\initial_< (I)$, \\ $<$ is a term ordering} $\Biggr\}$ 
$\leftrightarrow $
$\Biggl\{$ \parbox[c][2.5em][c]{4cm}{\centering marked Gröbner bases \\ of $I$} $\Biggr\}$ 
$\leftrightarrow$ 
$\Biggl\{$ \parbox[c][2.5em][c]{3.8cm}{\centering full-dimensional \\ cones of $\mathbb{G}(I)$} $\Biggr\}$
\end{center}
\end{theorem}
The first correspondence in \Cref{maintheorem} is immediate, whilst the second correspondence is a consequence of 
\cite[Theorem 1.11]{Sturmfels:1996}: a marked Gröbner basis $G$ encodes a facet description of the corresponding cone. 
Each element \((g,x^\alpha)\in G\) prescribes a leading term of \(g\) with respect to \(<\) which translate
to linear integral inequalities for the weight vectors \(\w\) given by
\begin{align*}
    \alpha\w \geq \beta\w, \quad \text{ where } x^\alpha = \inless(g) \text{ and } \beta \in \supp(g) \setminus\{\alpha\} \ \text{for every } g \in G,
\end{align*} where \(\supp(g)\) denotes the set of exponent vectors \(\alpha\in\mathbb{Z}^n\) such 
that \(x^\alpha\) has non-zero coefficient in \(g\).

Lower-dimensional cones in $\mathbb{G}(I)$ correspond to generalized initial ideals $\inw(I)$, where $\w$ is any weight vector in the relative interior of said cone. Generically, such ideals are \enquote{almost monomial}, and may be retrieved with the help of the following lemma: 
\begin{lemma}
    Let $G$ be a marked Gröbner basis of $I$ with respect to $<$ and $\w \in \Q^n _{\geq 0}$ be a weight vector on the boundary of the corresponding cone in $\mathbb{G}(I)$. The set 
    \[ \inw(G) = \{ \inw(g) , g \in G \} \]
    is a marked Gröbner basis of $\inw(I)$ with respect to $<$, where the marking of $\inw(g)$ is taken to be that of $g$.
\end{lemma}
At every step of the Gröbner walk, a basis of this form converted with Buchberger's algorithm and then lifted 
to the basis of $I$ corresponding to the adjacent full-dimensional cone, which corresponds to the refinement $<'_\w$. 
Recall that the \emph{normal form} of a non-zero polynomial $f \in R$ with respect to a Gröbner basis $G$ is the unique remainder obtained upon dividing $f$ by $G$ with the multivariate division algorithm (cf. \cite{Cox.Little.OShea:2015}, pg. 83). We denote this by $\overline{f}^G$.
\begin{lemma} \label{lem:lift}
    Let \(H\) be the Gröbner basis of \(I\) with respect to a term ordering \(<\) and let \(\w\in\QQ^n\)
    lie on the boundary of the cone in \(\mathbb{G}(I)\) corresponding to \(<\).
    
    If $M = \{ m_1 , ..., m_r \}$ is the marked Gröbner basis of $\inw(I)$ with respect 
    to the refinement ordering $<'_\w$, then
    \[ G := \{ m_1 - \overline{m_1}^{H} , ..., m_r - \overline{m_r}^{H} \} \] is a Gröbner basis of $I$ 
    with respect to $<'_\w$.
\end{lemma}
The Gröbner basis $G$ that is obtained from $M$ is referred to as the \emph{lift} of $M$ by $<$; 
it is not a marked Gröbner basis since it is not necessarily reduced.

In the setting of the Gröbner walk, the ordering $<$ from \Cref{lem:lift} is taken to be the start ordering, with respect to which a marked Gröbner basis is already known. The marked Gröbner basis which can be computed by reducing the lift Gröbner basis \(G\) corresponds to a full-dimensional cone in $\mathbb{G}(I)$ adjacent to the cone of \(<\) such that $\omega$ lies in their common intersection. As $\omega$ is chosen to be a point on the line segment from the cone of $<$ to that of the target ordering, the cone corresponding to $<'_\w$ will be \enquote{closer} to the cone of the target ordering in the Gröbner fan.  Thus, the process of subsequent passing to the generalized initial ideal and lifting to the adjacent basis may be repeated finitely many times to obtain a basis with respect to the target ordering. 
\section{Functionality}

Our implementation of the Gröbner walk ships with \texttt{OSCAR} since version~1.2.0, thus it
suffices to load \texttt{OSCAR}. There is a straightforward interface through the function \verb|groebner_walk|. 
\begin{example}\label{ex:running-example-gwalk}
    Continuing from example \Cref{ex:running-example}, we can calculate a Gröbner basis of the ideal
    \[
        I = \langle y^4+ x^3-x^2+x,x^4\rangle \triangleleft \QQ[x,y]
    \] with respect to \(\lexord\) by starting from a Gröbner basis for the \emph{graded reverse lexicographic ordering}
    \(\grevlexord\). In \texttt{OSCAR}, every polynomial ring comes with an internal ordering that is used
    for computations involving orderings. By default, this internal ordering is \(\grevlexord\).
    We adopt this interface and take the internal ordering as default if no other start ordering is specified.

    Thus, it suffices to call the Gröbner walk in the following way to calculate a Gröbner basis for 
    \(\lexord\) on \(\QQ[x,y]\) from a Gröbner basis with respect to \(\grevlexord\).
\begin{minted}{jlcon}
julia> groebner_walk(I, lex(R))
Gröbner basis with elements
  1: x + y^12 - y^8 + y^4
  2: y^16
with respect to the ordering
  lex([x, y])
\end{minted}
The corresponding computation path in the Gröbner fan is shown in \cref{fig:gfan-toy-example}.
\begin{figure}
    \centering
    \includegraphics[width=0.7\linewidth]{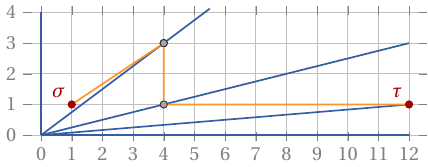}
    \caption{The Gröbner fan of the ideal \(\langle y^4+ x^3-x^2+x,x^4\rangle\). Each point 
    denotes an intermediate weight vector for which a Gröbner basis is computed. For implementation reasons, we choose an integer weight vector 
    in every intermediate step.}
    \label{fig:gfan-toy-example}
\end{figure}
\end{example}

When calling \verb|groebner_walk|, we construct the necessary initial information for the algorithm 
in the background and start the computation. This way, the user faces a single, 
uniform interface to each implemented algorithm. The selection of the algorithm used
is done by an (optional) \emph{keyword argument} \verb|algorithm=|.

Our implementation also offers additional diagnostic output to trace the computation. 
For example this allows us to inspect the intermediate steps of the computation in \Cref{ex:running-example-gwalk} 
which are shown in \cref{fig:gfan-toy-example}
\begin{minted}{jlcon}
julia> set_verbosity_level(:groebner_walk, 1);
julia> groebner_walk(I, lex(R))
Results for standard_walk
Crossed Cones in: 
ZZRingElem[1, 1]
ZZRingElem[4, 3]
ZZRingElem[4, 1]
ZZRingElem[12, 1]
Cones crossed: 4
Gröbner basis with elements
  1: x + y^12 - y^8 + y^4
  2: y^16
with respect to the ordering
  lex([x, y])
\end{minted}

The default choice is the variant by \citeauthorandcite{Collart:GWalk:1997} (which can also 
be specified by setting \verb|algorithm=:standard|). Another supported choice is 
\verb|:generic| for the \emph{generic walk} by \citeauthorandcite{GenericWalk:2007}. 

The default choices for starting and target ordering are the internal ordering of 
\texttt{R} and the lexicographical ordering respectively. If not specified further,
this internal ordering is usually the \emph{degree reverse lexicographical ordering}.
In general, it is not obvious though what is a sensible choice for a starting ordering.
One possible way to alleviate this is to use \emph{Gröbner basis detection methods},
which we discuss in \Cref{sec:outlook}.

As a word of caution, since we rely on the implementation of Gröbner bases in \texttt{OSCAR} which in turn 
is based on the facilities for Gröbner bases in Singular, we also inherit the restrictions on computations with
Gröbner bases. Namely, the backend in Singular only supports weights less than \(2^{31}-1\).
This means that Gröbner walk computations which involve a significant blowing up of weight vectors 
throw an \texttt{InexactError}. At present, we do not handle this limitation.

\section{Technical Contribution}
We implement two variants of Gröbner walk algorithms, the standard walk by 
\citeauthorandcite{Collart:GWalk:1997} and the generic walk by \citeauthorandcite{GenericWalk:2007}.
For the generic walk, we also provided a naive implementation of marked Gröbner bases.

Problems of Gröbner basis conversion arise in a wide variety of contexts. 
By default, the \texttt{OSCAR} \cite{OSCAR} function \verb|groebner_basis| computes a Gröbner basis using Buchberger's algorithm.  
This approach did not terminate in reasonable time in our larger examples.

Several other Gröbner basis algorithms have been implemented in \texttt{OSCAR} 
\cite{Faugere:1993,Faugere:1999,Simoes:2014}, and may be used via the keyword \verb|algorithm=| when calling \verb|groebner_basis|.
While all of these methods are improvements on Buchberger's algorithm, they each come with their own limitations; 
for example, FGLM \cite{Faugere:1993} is only applicable to zero-dimensional ideals, 
whereas the current \texttt{OSCAR} implementation of the F4 algorithm \cite{Faugere:1999} may be called for ideals
over the rationals or finite fields of machine-size characteristic, but only calculates Gröbner bases for
\(\grevlexord\).
In contrast to this, the Gröbner walk (and our implementation in \texttt{OSCAR}) works in full generality; it may be called on ideals over $\Q$ or $\mathbb{F}_p$ of arbitrary dimension and for arbitrary term orderings. 

This last fact makes it especially well-suited for problems of elimination. 
Furthermore, due to the variable and unpredictable performance of Gröbner basis computations on generic ideals, it is advantageous to have a variety of options for these tasks. 

Our implementation is included in version~1.5.0 of \texttt{OSCAR}  
as an \emph{experimental package}, which means our implementation is be shipped 
with \texttt{OSCAR} as a submodule.
Alternatively, the implementation is provided as self-contained Julia package with dependency on \texttt{OSCAR}
which can be found at \cite{GroebnerWalk.jl}. Calling the function \texttt{GroebnerWalk.groebner\_walk} calls the GroebnerWalk.jl version of our implementation, which may differ from the OSCAR version in the future.

\section{Comparison to classical Gröbner basis algorithms}
For our comparisons we choose two types of problems and ran computations over \(\QQ\) and \(\mathbb{F}_p\)
for \(p = 11863279\).
The first kind are computations of lexicographic Gröbner bases of zero-dimensional ideals for solving systems of polynomial equations. The chosen systems are from Jan Verschelde's database 
\footnote{The database can be found at \url{http://homepages.math.uic.edu/~jan/demo.html}.}
and included in the software \verb|PHCpack| \cite{Verschelde:1999}. 
They are commonly used in the benchmarks of polynomial solvers.
The second kind are computations of Gröbner bases of ideals of dimension $\geq 1$ with respect to \emph{elimination term orderings}. These computations arise in problems of \emph{implicitization} of surfaces given in parametric form, such as in the \texttt{agk} \cite{Amrhein:1997} and \texttt{newell} \cite{Tran:2004a} examples. 

\begin{table}
    \caption{Overview of the polynomial systems \(I = \langle G\rangle\) chosen for the comparison
    with a summary of characteristics.}
    \label{tab:test-candidates}
    \begin{tabular}{llccccccc}
        \toprule
        \multicolumn{3}{c}{Polynomial systems} & \multicolumn{5}{c}{Characteristics} \ \\\cmidrule(r){1-3}\cmidrule(l){4-8}
        Name & Description & Ref & \(<_\sigma\) & \(<_\tau\) & \(\lvert G\rvert\) & \(\lvert G_\sigma\rvert\) & \(\lvert G_\tau\rvert\)\\\midrule
        \texttt{cyclic5} & The cyclic 5-roots problem & \cite{Bjorck:1991} & degrevlex & lex & 5 & 20 & 30  \\
        \texttt{cyclic6} & The cyclic 6-roots problem &  &  &  & 6 & 45 & 70  \\\midrule
        \texttt{katsura6} & A problem of magnetism in physics & \cite{Katsura:1990} & degrevlex & lex & 7 & 41 & 64 \\
        \texttt{katsura7} &                                   &           &  &  & 8 & 74 & 128 \\
        \texttt{katsura8} &                                   &           &  &  & 9 & 143 & 256 \\\midrule
        \texttt{agk4} & A parametric B\'ezier surface & \cite{Amrhein:1997} & \eqref{elimination-orders} & \eqref{elimination-orders} & 3 & 3 & 29 \\
        \texttt{newell} & The Newell teapot & \cite{Tran:2004a} & \eqref{elimination-orders} & \eqref{elimination-orders} &  3 & 12 & 39 \\
        \texttt{tran3.3} & Example 3.3 from \citeauthor{Tran:2000} & \cite{Tran:2000} & degrevlex & lex & 2 & 5 & 10
        \\\bottomrule
    \end{tabular}
\end{table}

Following \cite{Tran:2004a} we chose the start and target orderings for those two problems as represented by the matrices 
\begin{equation}\label{elimination-orders}
    {<_\sigma} = \begin{pmatrix}
      1 & 1 & 1 & 0 & 0 \\
      0 & 0 & 0 & 1 & 1 \\
      0 & 0 & 0 & 1 & 0 \\
      1 & 1 & 0 & 0 & 0 \\
      1 & 0 & 0 & 0 & 0
      \end{pmatrix}
    \quad\text{and}\quad
    {<_\tau} = \begin{pmatrix}
      0 & 0 & 0 & 1 & 1  \\
      1 & 1 & 1 & 0 & 0 \\
      1 & 1 & 0 & 0 & 0 \\
      1 & 0 & 0 & 0 & 0 \\
      0 & 0 & 0 & 1 & 0
      \end{pmatrix}.
\end{equation}

We ran the comparisons on a MacBook Pro with an \num{2,4}\unit{\giga\hertz} Apple M2 Max. Each computation was allotted a maximum of
\num{16}\unit{\giga\byte} of memory.
We used macOS 26.0 with Julia 1.10.10 and \texttt{OSCAR} 1.5.0. The results of the 
comparison are shown in \Cref{tab:timings}. The data and code to run the benchmarks 
can be found at \url{https://zenodo.org/records/17473425}.

Owing to the high upper bounds in Gröbner basis computations and the resulting unpredictability, we report mixed results 
in our comparisons. While the implementation of classical Gröbner basis algorithms in
\texttt{OSCAR} performs reasonably well on the zero-dimensional \texttt{cyclic} ideals and some instances from
implicitization, there are also examples where the Gröbner walk performed better or produced a Gröbner basis at all.
We believe that the poor performance of the generic walk despite its theoretical advantages 
is due to the currently sub-optimal implementation. For example, the reduction steps are still performed using naive 
polynomial long division without any of the possible algorithmic improvements. 
Comparison with the Macaulay2 implementation \cite{M2Gwalk} yielded mixed results, despite the implementations being very similar on a theoretical level. The results are presented in \Cref{tab:M2comparisons}. Interestingly, none of the computations involving the \texttt{katsura} polynomial systems terminated in Macaulay2. We attribute the generally superior performance of the \texttt{OSCAR} implementation of the standard walk to the more optimized methods for polyhedral geometry and Gröbner basis computation via the Polymake and Singular backends respectively. The fact that the Macaulay2 implementation of the generic walk performs better than \texttt{OSCAR} in certain instances is likely due to the native \texttt{markedGB} data structure, which allows for the encoding of marked Gröbner bases in Macaulay2 without the explicit specification of a term ordering.
\begin{table}
    \caption{Average runtimes for the Gröbner basis computations in Table 1 using our implementations of the Gröbner walk as well as OSCAR's built-in \texttt{groebner\_basis} function. Missing entries indicate a computation that
    timed out. The cutoff was \num{3000}\unit\second. \(p = 11863279\).}
    \begin{tabular}{lrrrrrr}
        \toprule
        &\multicolumn{6}{c}{Runtime}\\\cmidrule{2-7}
        System & \multicolumn{2}{c}{Standard walk} & \multicolumn{2}{c}{Generic walk} &  \multicolumn{2}{c}{\texttt{groebner\_basis}}\\\cmidrule(lr){2-3} \cmidrule(lr){4-5} \cmidrule(lr){6-7}
        & \multicolumn{1}{c}{\(\QQ\)} & \multicolumn{1}{c}{\(\mathbb{F}_p\)} & \multicolumn{1}{c}{\(\QQ\)} & \multicolumn{1}{c}{\(\mathbb{F}_p\)} & \multicolumn{1}{c}{\(\QQ\)} & \multicolumn{1}{c}{\(\mathbb{F}_p\)} \\\midrule
        \texttt{cyclic5} & 
        \num{8.35}\unit{\milli\second} & \num{5.29}\unit{\milli\second} & 
        \num{0.13}\unit\second & \num{0.23}\unit\second & 
        \num{3.82}\unit{\milli\second} & \num{0.68}\unit{\milli\second} \\
        \texttt{cyclic6} & 
        \num{0.17}\unit\second & \num{0.04}\unit\second & 
        \num{3.42}\unit\second & \num{7.67}\unit\second & 
        \num{0.08}\unit\second & \num{0.04}\unit\second \\\midrule
        \texttt{katsura6} & 
        \num{0.20}\unit\second & \num{0.06}\unit\second & 
        \num{13.23}\unit\second & \num{10.93}\unit\second & 
        \num{20.72}\unit\second & - \\
        \texttt{katsura7} & 
        \num{2.37}\unit\second & \num{0.58}\unit\second & 
        \num{365.80}\unit\second & \num{156.02}\unit\second & 
        \num{59.80}\unit\second & - \\
        \texttt{katsura8} & 
        \num{23.35}\unit\second & \num{4.92}\unit\second & 
        - & - & 
        -  & - \\\midrule
        \texttt{agk4} & 
        \num{3.98}\unit\second & \num{0.71}\unit\second & 
        \num{39.79}\unit\second & \num{23.95}\unit\second & 
        \num{0.55}\unit\second & \num{0.10}\unit\second \\
        \texttt{newell} & 
        \num{20.42}\unit\second & \num{3.17}\unit\second & 
        - & \num{833.53}\unit\second & 
        \num{808.14}\unit\second & \num{0.16}\unit\second \\
        \texttt{tran3.3} & 
        \num{0.40}\unit\second & \num{0.18}\unit\second & 
        \num{3.40}\unit\second & \num{0.66}\unit\second & 
        \num{0.05}\unit\second & \num{3.35}\unit{\milli\second}
        \\\bottomrule
    \end{tabular}
    \label{tab:timings}
\end{table}

\begin{table} 
    \centering
    \caption[Side-by-side comparison of the Gröbner walk performance in \texttt{Macaulay2} and \texttt{OSCAR}]{ Side-by-side comparison of the standard and generic walks in \texttt{OSCAR} and \texttt{Macaulay2} over $\Q$ and $\mathbb{F}_p$ ($p = 11863279$, cutoff 3000s) }
    \begin{tabular}{lrrrrrr}
        \toprule
        &\multicolumn{6}{c}{avg. runtime  over $\Q$}\\\cmidrule{2-7}
         System & \multicolumn{2}{c}{Standard walk} & \multicolumn{2}{c}{Generic walk} &  \multicolumn{2}{c}{\texttt{groebner\_basis} / \texttt{gb}}\\
        \cmidrule(lr){2-3} \cmidrule(lr){4-5} \cmidrule(lr){6-7}
        & \texttt{OSCAR} & \texttt{M2} & \texttt{OSCAR} & \texttt{M2} & \texttt{OSCAR} & \texttt{M2} \\
        \texttt{cyclic5} &    \num{8.35}\unit{\milli\second} & \num{0.02}\unit{\second} & 
        \num{0.13}\unit\second & \num{0.25}\unit\second& 
        \num{3.82}\unit{\milli\second} & \num{16.2}\unit{\milli\second} \\
        
        \texttt{cyclic6} & \num{0.17}\unit\second & \num{0.69}\unit\second & 
        \num{3.42}\unit\second & \num{4.29}\unit\second & 
        \num{0.08}\unit\second & \num{0.55}\unit\second  \\
        \cmidrule{1-7}
        \texttt{katsura6} & \num{0.20}\unit\second & - & 
        \num{13.23}\unit\second & - & 
        \num{20.72}\unit\second & - \\
        \cmidrule{1-7}
        \texttt{agk4} & \num{3.98}\unit\second & \num{0.93} \unit\second& 
        \num{39.79}\unit\second & \num{2.19}\unit\second& 
        \num{0.55}\unit\second & \num{0.75}\unit\second   \\
        \cmidrule{1-7}
        &\multicolumn{6}{c}{avg. runtime (s.) over $\mathbb{F}_p$}\\\cmidrule{2-7}
        \texttt{cyclic5} & \num{5.29}\unit{\milli\second} & \num{0.27}\unit\second  & \num{0.23}\unit\second & 
          \num{4.16}\unit\second & \num{0.68}\unit{\milli\second} & \num{3.23}\unit{\milli\second} \\ 
        
        \texttt{cyclic6} & \num{0.04}\unit\second & \num{0.27}\unit\second & 
          \num{7.67}\unit\second & \num{4.19}\unit\second &
          \num{0.04}\unit\second & \num{0.12}\unit\second  \\
        \cmidrule{1-7}
        \texttt{katsura6}  & \num{0.06}\unit\second & - & 
         \num{10.93}\unit\second & - & 
        - & - \\
        \cmidrule{1-7}
        \texttt{agk4} & \num{0.71}\unit\second & \num{2.97}\unit\second & 
         \num{23.95}\unit\second & \num{1.69}\unit\second &
         \num{0.10}\unit\second & \num{0.12}\unit\second   \\
        \bottomrule
    \label{tab:M2comparisons}
    \end{tabular}
\end{table}

\section{Future directions}\label{sec:outlook}
To compare the performance of the standard walk to the generic walk, we analyzed 
profiling data for the \texttt{tran3.3} and \texttt{agk4} instances. 
Those examples have a comparably big difference in computation time for standard and generic
walk while the runtime of the generic walk is still low enough to reasonably compute 
traces.

In the generic walk computations the most significant bottleneck is the calculation of normal forms with respect to the symbolic
intermediate orderings. This is due to the fact that we had to implement division with remainder naively
as the existing implementations require at least a weight ordering with known weight vector.
The next step would be to adopt a linear algebra approach to computing normal forms using Macaulay matrices.
However, there is currently no user-facing \texttt{OSCAR} function for this and in any case, 
such a method would have to be adapted to be compatible with our marked Gröbner basis structure.  
Also, both walks ultimately rely on the Singular for Gröbner basis computations, which at the time of writing does 
not support computations with weight orderings with entries larger than $2^{31} -1$. 
Currently, the only work-around would be to write new Gröbner basis functions in \texttt{OSCAR} 
(which itself supports arbitrary precision integers) without any interaction with Singular, 
the implementation of which is a worthwhile avenue itself to extend the functionality of \texttt{OSCAR} but lies beyond our scope. 

An important point in any variant of the Gröbner walk is the choice of starting term ordering
to compute a first Gröbner basis. While rare, the specified generators of an ideal might 
already form a Gröbner basis with respect to some term ordering. 
There exist criteria for \emph{Gröbner basis detection}, see \eg Chapter 3 in \cite{Sturmfels:1996}. 
By employing such criteria, we can avoid the first computation of
a Gröbner basis. A detection criterion has been implemented in Julia by 
\citeauthorandcite{Borovik.Duff.ea:2024}.

\section*{Acknowledgments}
We would like to thank Carlos Améndola, Anne Frühbis-Krüger, Ben Hollering, Michael Joswig, Yue Ren, as well as an anonymous referee
for helpful comments, guidance and pointers on how to improve our implementation.

KF is funded by the Deutsche Forschungsgemeinschaft (DFG, German Research Foundation) 
under Germany´s Excellence Strategy – The Berlin Mathematics Research Center MATH+ 
(EXC-2046/1, project ID: 390685689).
FN was supported by the SPP 2458 \enquote{Combinatorial Synergies}, funded by the Deutsche Forschungsgemeinschaft (DFG, German Research Foundation) – project ID: 539875257. 

\singlespacing
\printbibliography

\end{document}